\documentclass[twocolumn]{svjour3}         
\usepackage{graphicx}
\usepackage{amsmath}
\usepackage{shortcuts}
\usepackage{amssymb}
\usepackage{epsf}
\usepackage{cuted}
\usepackage{subeqn}

\journalname{Nonlinear Dynamics}
\begin{document}

\title{Stability of a 2-dimensional Mathieu-type system with quasiperiodic coefficients}

\author{Thomas J. Waters}

\institute{Thomas J. Waters \at
              Dept. Applied Mathematics \\
              NUI Galway, Ireland
              \email{thomas.waters@nuigalway.ie}
}

\date{Received: date / Accepted: date}

\maketitle

\begin{abstract}
In the following we consider a 2-dimensional system of ODE's containing quasiperiodic terms. The system is proposed as an extension of Mathieu-type equations to higher dimensions, with
emphasis on how resonance between the internal frequencies leads
to a loss of stability. The 2-d system has two `natural'
frequencies when the time dependent terms are switched off, and it
is internally driven by quasiperiodic terms in the same
frequencies. Stability charts in the parameter space are generated
first using numerical simulations and Floquet theory. While some
instability regions are easy to anticipate, there are some
surprises: within instability zones small islands of stability
develop, and unusual `arcs' of instability arise also. The
transition curves are analyzed using the method of harmonic
balance, and we find we can use this method to easily predict the `resonance curves'
from which bands of instability emanate. In addition, the method of multiple scales is used to examine the
islands of stability near the 1:1 resonance. 
\end{abstract}

\section{Introduction}

The study of systems of ode's with periodic coefficients arises naturally in many branches of applied mathematics.
Aside from being interesting in their own right, in applications such systems arise when considering the orbital
stability of a periodic solution to some dynamical system (see, for example, Betounes \cite{betounes}). The archetypal equation for this class is the Mathieu
equation, \be \ddot{x}+(\delta+\vep\cos t)x=0 \label{math1d} \ee where a dot denotes differentiation w.r.t. time.
This equation, though linear, is sufficiently complicated by its explicit time dependence that closed form solutions
do not exist. Instead, the common method of attack is to
consider the stability of the origin using Floquet theory (or other methods). The ($\delta,\vep$) parameter space
is clearly divided into regions of stability/instability (see Jordan and Smith \cite{jordan} for a detailed discussion).
The appearance of these instability zones is intuitively understood in the following way: consider the time-dependent
term in \eqref{math1d} as a perturbation. When $\vep=0$, $x$ has periodic solutions with frequency $\sqrt{\delta}$
(the `natural' frequency of the system). As we switch on the perturbation, there will be a resonance between the
`driving' frequency, 1, and the natural frequency $\sqrt{\delta}$, which leads to instability.

There are many extensions to this equation, of which we describe
only a few. Firstly the equation can be extended to many
dimensions, by replacing $x$ with a vector and $\delta$ and $\vep$
with matrices. This problem was analyzed in Hansen \cite{hansen}
using Floquet theory, however (in contrast with the present work)
the time dependent terms only depended on one frequency. Secondly
the equation may be extended to contain nonlinear terms; for
example El-Dib \cite{eldib} consider an extension up to cubic
order with each term containing subharmonics in the periodic term,
and Younesian et al \cite{youn} consider also cubic order with two
small parameters and terms in $\sin$ and $\cos$.

Thirdly, and most importantly for the present work, we may develop
the time dependent term to contain two periodic terms in different
frequencies, and as the two frequencies need not be commensurable
this system is referred to as the `quasiperiodic Mathieu
equation'. This type of system was developed and studied in great
detail by Rand and coworkers: in Rand and Zounes \cite{rand1} and
Rand et al \cite{rand4} the following quasiperiodic Mathieu
equation is examined \be \ddot{x}+[\delta+\vep(\cos(t)+\cos(\omega
t))]x=0 \label{quasi1d} \ee and transition curves in the
$(\delta,\omega)$ parameter space are developed. Specific
resonances are examined in Rand and Morrison \cite{rand2} and Rand
et al \cite{rand3}. In Sah et al \cite{rand5} a form of
the quasiperiodic Mathieu equation is derived and examined in the
context of the stability of motion of a particle constrained by a
set of linear springs, and finally in Rand and Zounes \cite{rand6} a
nonlinear version of the quasiperiodic Mathieu equation is studied. In these papers a combination of numerical
and analytical techniques are used to examine how resonances
between the internal driving frequencies lead to a loss of
stability in solutions, and it is very much in the spirit of these
papers that we present the following work.

We seek to extend consideration of quasiperiodic Mathieu-type
systems to higher dimensions. The motivation for this is clear: as
mentioned previously, linear systems with (quasi)periodic
coefficients often arise when a nonlinear dynamical system is
linearized around a (quasi)periodic solution, the so-called
`variational equations' \cite{betounes}. As most dynamical systems
of interest will have more than one degree of freedom (for example
the motion of a particle in space under the influence of a
collection of forces), the associated variational equations will
also have more than one degree of freedom (unless only certain
modes are being considered as in \cite{rand5}). We seek a
2-dimensional system containing more than one frequency and a
small parameter $\vep$ to enable a perturbation analysis. However,
it is beneficial to keep the number of free parameters low in
order to make the behavior clear. With this in mind, we shall
consider the following system: \be
\ddot{x}+[\alpha^2+\vep\cos(\beta t)]x+\vep\cos(\alpha t)y =0
\nonumber \\ \ddot{y}+\vep \cos(\beta t)x+[\beta^2+\vep
\cos(\alpha t)]y=0 \label{syseqn} \ee which in general is of the
form \be \frac{d^2 X}{dt^2}+[P_0+\vep P_1(t)]X=0 \label{matform}
\ee with $X=(x,y)$ and \be P_0=\left(
                               \begin{array}{cc}
                                 \alpha^2 & 0 \\
                                 0 & \beta^2 \\
                               \end{array}
                             \right), \quad P_1(t)=\left(
                               \begin{array}{cc}
                                 \cos(\beta t) & \cos(\alpha t) \\
                                 \cos(\beta t) & \cos(\alpha t) \\
                               \end{array}
                             \right). \ee Here a dot denotes
                             differentiation w.r.t. time, $\vep$
is small and $\alpha,\beta \in \mathbb{R}$.

We can predict to some degree how solutions of this system will
behave: when $\vep=0$, the solutions are periodic in the natural
frequencies $\alpha$ and $\beta$. As the perturbation, $\vep
P_1(t)X$, is switched on we would expect resonance between the
natural and driving frequencies, as in the Mathieu equation
\eqref{math1d}. As the natural and driving frequencies are the
same, this means we expect instability to arise where $\alpha$ and
$\beta$ are commensurable. And this is to a degree what we do see,
with two unexpected features. Firstly, curved `arcs' of
instability can be seen which clearly do not arise from integer
ratios of the parameters. Secondly, within the bands of
instability small islands of stability appear. These two features
will be addressed in more detail in what follows.

While this system may seem a little artificial, we will show how
it may arise as an approximate variational equation. Consider the
following system of nonlinear ODE's: \be \ddot{u}=-\alpha^2
u-uv,\quad \ddot{v}=-\beta^2 v-uv. \label{auton} \ee This is an example of a system which is not dissipative but nonetheless is not Hamiltonian; as such it falls into the more general class of Louiville systems (see below). The author's direct experience with systems of this type is in the problem of the solar sail in the circular restricted 3-body problem (see for example the author and McInnes \cite{solar}), although other applications exist.

The origin of \eqref{auton} is a fixed point, and letting $u=\tilde{u}+x,\ v=\tilde{v}+y$ we have the
linear system \be \frac{d^2X}{dt^2}=
      \left(
        \begin{array}{cc}
          -\alpha^2-v & -u \\
          -v & -\beta^2-u \\
        \end{array}
      \right)_{u=\tilde{u},v=\tilde{v}}\ X \label{autonlin} \ee where $X=(x,y)$. Letting $(\tilde{u},\tilde{v})=(0,0)$
      we find two linear oscillators which we write as $\vep\cos(\alpha t),\
      \vep\cos(\beta t)$ (letting the amplitudes be equal and setting the phases to zero),
      and then letting $\tilde{u}=\vep\cos(\alpha t),\ \tilde{v}=\vep\cos(\beta t)$ we
      get system \eqref{syseqn} above. Note that the autonomous system \eqref{auton} is not
      Hamiltonian and therefore \eqref{autonlin} is not symmetric;
      neither, for the same reason, is $P_1(t)$ in \eqref{syseqn}.
      Whether or not the system under consideration is Hamiltonian
      only becomes an issue in dimension 2 and above; in one
      dimension we can always find a (time dependent) Hamiltonian,
      a fact taken advantage of in Ye\c{s}ilta\c{s} and \c{S}im\c{s}ek \cite{yessim}.

Note also that $(u,v)=(\vep\cos(\alpha t),\vep\cos(\beta t))$ is
{\it not} a solution to \eqref{auton}, rather a solution to the
linearized system \eqref{autonlin}. As such \eqref{syseqn} is not
a true variational equation, rather an approximation to one. We
make this approximation, rather than seek more accurate solutions
to \eqref{auton}, so as to enable an analytical rather than purely
numerical treatment. The Mathieu equation itself can be seen in
this light, as an approximate variational equation to the system
\[ \ddot{u}=-u-u^3. \] The fact that \eqref{syseqn} is only
approximately a variational equation will be of relevance in what
follows.

The rest of the paper is laid out as follows: In Section 2 we
numerically integrate system \eqref{syseqn} to examine the
stability properties of various values of $\alpha$ and $\beta$. As
computers can only handle finite precision we must choose rational
values of $\alpha$ and $\beta$ and therefore \eqref{syseqn} is
periodic, rather than quasiperiodic (a more appropriate term might
be multiply-periodic). As the system is periodic, albeit of
possibly very long period, we may use Floquet theory to enhance
the numerical integration. We are able to generate stability
charts in the $(\alpha,\beta)$ parameter space for various values
of $\vep$ and we discuss the main features of same. The fact that
we have lost the true quasiperiodic nature of the system in this
section is made up for by the following: every irrational number
has rational numbers arbitrarily close to it which thus give a
very good approximation; the numerically generated stability
charts allow us to isolate the interesting features of the system
which informs later analysis; and the remaining analysis in the
paper holds for arbitrary $\alpha$ and $\beta$ and thus recaptures
the quasiperiodic nature of the system.

In Section 3 we use the method of harmonic balance to find the
transition curves in parameter space. Writing solutions in terms of truncated Fourier series we approximate the transition curves by the vanishing of the Hill's determinant. As the truncation size becomes large the matrices of vanishing determinant also become large making this method cumbersome.
Also, there are convergence issues with the infinite Hill's
determinant appropriate to \eqref{syseqn} which we will elaborate
on. However, we show that we can use the method of harmonic balance to quickly find approximations to transition curves by understanding the following mechanism of instability: when $\vep=0$ there are certain `resonance curves' in parameter space (for the Mathieu equation \eqref{math1d} these are simply the points $\delta=n^2/4$ for $n\in\mathbb{Z}$). As $\vep$ grows, these curves widen into bands of
instability. The method of harmonic balance allows us to find
these resonance curves quickly, and we demonstrate the usefulness
of this technique on a second system closely related to
\eqref{syseqn}.

In Section 4 we use the method of multiple scales to examine in
detail the fine structure near the $\alpha\approx\beta$ resonance.
From Section 1 we can see small islands of stability arising in
the middle of a large band of instability, and we examine this
feature more closely. We separate out the `slow flow' using three
time variables, and find we can analytically predict the
appearance of these pockets of stability. Finally we give some
discussion and suggestions for future work.

\section{Numerical integration}

We wish to fix the values of the parameters $(\alpha,\beta,\vep)$
and numerically integrate the system \eqref{syseqn} to examine the
stability of the origin. In deciding what range of values to
choose for the parameters, we note that \eqref{syseqn} has the
following scale invariance: consider the system \be
\ddot{x}+[(m\alpha)^2+m^2\vep\cos(m\beta t)]x+m^2\vep\cos(m\alpha
t)y =0 \nonumber \\ \ddot{y}+m^2\vep \cos(m\beta
t)x+[(m\beta)^2+m^2\vep \cos(m\alpha t)]y=0  \ee By defining a new
time coordinate  $t'=mt$ we recover the original system
\eqref{syseqn}. Thus the following parameter values are equivalent
as regards stability: \[ (\alpha,\beta,\vep),\quad
(m\alpha,m\beta,m^2\vep). \] This means we need only consider
$\alpha$ and $\beta$ on the unit square (for example), and the
stability of other values of $\alpha$ and $\beta$ are found by
varying the value of $\vep$. For example, the stability of system
\eqref{syseqn} when $\alpha=3,\beta=2,\vep=0.1$ is the same as
when $\alpha=3/4,\beta=1/2,\vep=0.0125$.

As mentioned in the introduction, to integrate numerically we must
use finite precision values of $\alpha$ and $\beta$ and thus the
system is (multiply-)periodic. As such, we may take advantage of
Floquet theory.

\begin{center}
\begin{figure*}[htb!]
  \includegraphics[width=0.85\textwidth]{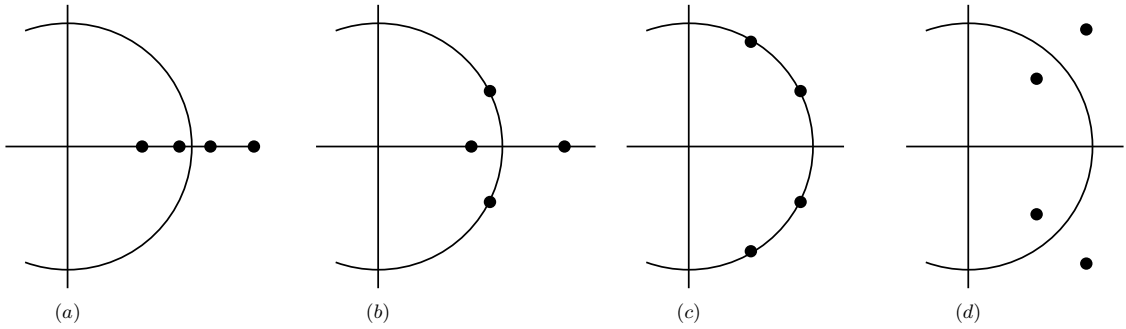}
\caption{Possible positions of multipliers in the complex plane with respect to the unit circle
for a $t$-invariant system of order 4.} \label{multipliers}
\end{figure*}
\end{center}

Writing \eqref{syseqn} as a first order system, the fundamental matrix $\Phi(t)$ solves the following initial value
problem (see, for example, Yakubovich and Starzhinskii \cite{yak}
for a detailed treatment of Floquet theory) \be \frac{d
\Phi(t)}{dt}=A(t)\Phi(t), \quad \Phi(0)=I \label{lininit} \ee
where \be A(t)= \left(
       \begin{array}{cccc}
         0 & 0 & 1 & 0 \\
         0 & 0 & 0 & 1 \\
         -[\alpha^2+\vep \cos(\beta t)] & -\vep\cos(\alpha t) & 0 & 0 \\
         -\vep\cos(\beta t) & -[\beta^2+\vep \cos(\alpha t)] & 0 & 0 \\
       \end{array}
     \right) \label{linsys} \ee and $A(t+T)=A(t)$ for some fundamental period $T$. The
     stability of the solutions to this system are determined by the monodromy matrix\footnote{In some of the literature, $A(t)$ is known as the monodromy matrix and $\Phi(T)$ as the matrizant, propagator or characteristic matrix.}, given by
     $\Phi(T)$, since it can be shown that \[ \left(
                                                      \begin{array}{c}
                                                        x(T) \\
                                                        y(T) \\
                                                      \end{array}
                                                    \right)
     =\Phi(T)\left(
                                                      \begin{array}{c}
                                                        x(0) \\
                                                        y(0) \\
                                                      \end{array}
                                                    \right).
     \] More precisely, stability is determined by the eigenvalues of $\Phi(T)$, the so-called
     (characteristic) multipliers $\lambda_i$. If $|\lambda_i|\leq
     1$ the solutions remain bounded, implying stability, and if
     any $|\lambda_i|>1$ the solution is unbounded.

From the structure of \eqref{linsys} we can make some comments
about the multipliers:
\begin{enumerate}
  \item Louiville's theorem \cite{jordan} gives
  \[ det[\Phi(T)]= \lambda_1 \lambda_2 \lambda_3 \lambda_4=
  \exp\left(\int_0^T Tr\{A(t)\}dt\right)=1. \label{louiv} \] This however
  does {\it not} mean that the eigenvalues must appear in reciprocal
  pairs, merely that all four multiplied together must be 1.
  \item That they must appear in reciprocal pairs comes from the {\it t-invariant}
  nature of the system \cite{yak}. These are systems with the property
  \be A(-t)G+GA(t)=0 \label{tinv} \ee for some non-singular constant matrix $G$.
  Since $A(t)$ is even, i.e.\ $A(-t)=A(t)$, this holds with \be G=\left(
                                                                        \begin{array}{cc}
                                                                          I & 0 \\
                                                                          0 & -I \\
                                                                        \end{array}
                                                                      \right), \ee
      the entries understood to be $2\times 2$ matrices. In this case, we can show
      \be \Phi(T)^{-1}=G. \Phi(T). G^{-1}, \ee in other words, the monodromy matrix
      and its inverse are similar and thus have the same eigenvalues; thus if
      $\lambda_i$ is a multiplier then so is $1/\lambda_i$.

      As $\Phi(T)$ is real, complex eigenvalues must appear in
      conjugate pairs, and since they must also be reciprocals
      this means complex eigenvalues must be on the unit circle
      (with one rare exception, see point 4 below). Therefore,
      complex multipliers represent stability, and real
         multipliers instability.
  \item If the system \eqref{syseqn} were a variational equation, that is a
  (nonlinear autonomous) system linearized around a periodic solution, we would
  expect a unit eigenvalue (see chap. 7 of Betounes \cite{betounes}) and therefore (from the previous point) 2 unit eigenvalues.
  This would make life easier, as stability would be ensured if \be Tr\{\Phi(T)\}\leq 4 \ee
  and thus we would not need to calculate the multipliers themselves. Unfortunately we
  cannot say this is the case (see the Introduction), and so there is no guarantee of unit eigenvalues. As the system can be viewed as an approximation
   to a variational equation then we would expect (for small $\vep$) to have two multipliers close to, but not equal to, one.
  \item There are a number of possibilities for the positions of the multipliers in
  the complex plane, shown in Figure \ref{multipliers}. The most common is for the set
  \be 1+e,\quad \frac{1}{1+e},\quad a+bi,\quad a-bi \ee with $e\in\mathbb{R}^+$
  and $a^2+b^2=1$ (since $\bar{\lambda}=1/\lambda$). We may describe this as `saddle-centre'
  (Figure \ref{multipliers}(b)). Stability is given by `centre-centre' configurations (Figure \ref{multipliers}(c)),
  however we must also accept the possibility that the multipliers may not be on the unit circle (Figure \ref{multipliers}(d)).
  This is possible when the eigenvalues of the centre-centre configuration
  come together and then bifurcate off the unit circle (a Krein collision). Note, the intermediate stage here means the
  eigenvalues are equal and this is a degenerate case; this will not be pursued further.
\end{enumerate}

\begin{figure*}
\begin{center}
\resizebox*{4.7in}{!}{\includegraphics{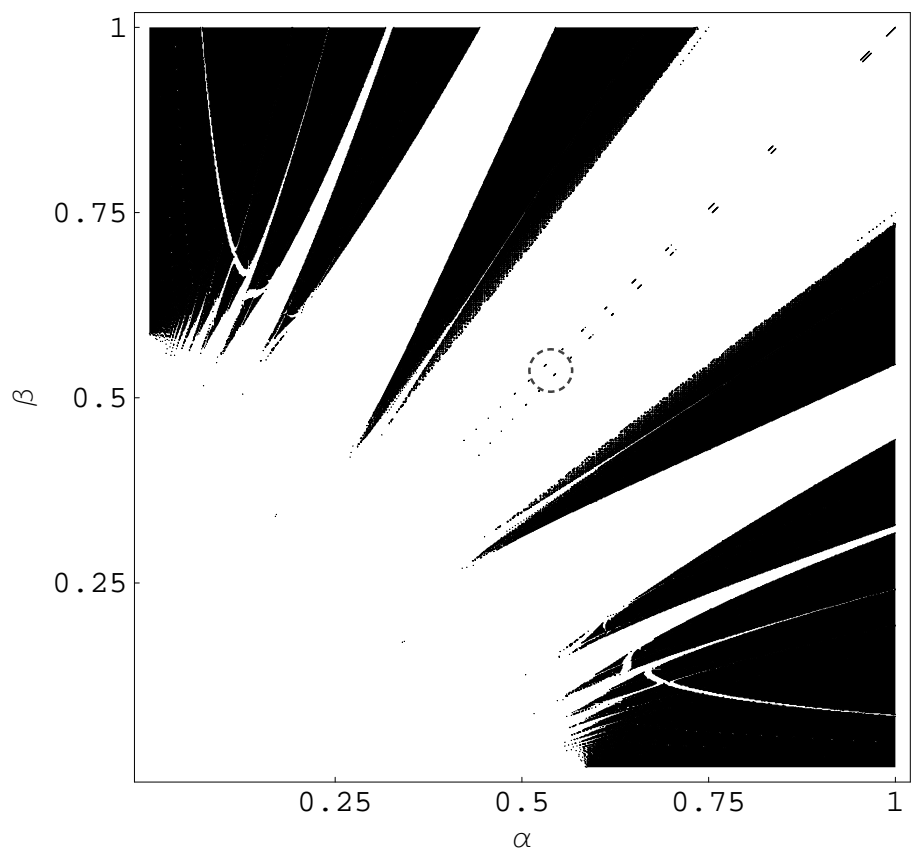}}
\end{center}
\caption{Stability diagram for system \eqref{syseqn} with
$\vep=0.1$ and the $\alpha$, $\beta$ unit square divided into an
$800 \times 800$ grid. White represents points where the norm of
the maximum multiplier is greater than 1.025 (and therefore the
system is unstable), and black where the norm is less than 1.025
(the system being stable). The dashed circle shows the region of
parameter space where the solutions in Figure \ref{dpoint} are
taken from.} \label{stabchart}
\end{figure*}

\begin{figure*}
\begin{center}
\resizebox*{3.5in}{!}{\includegraphics{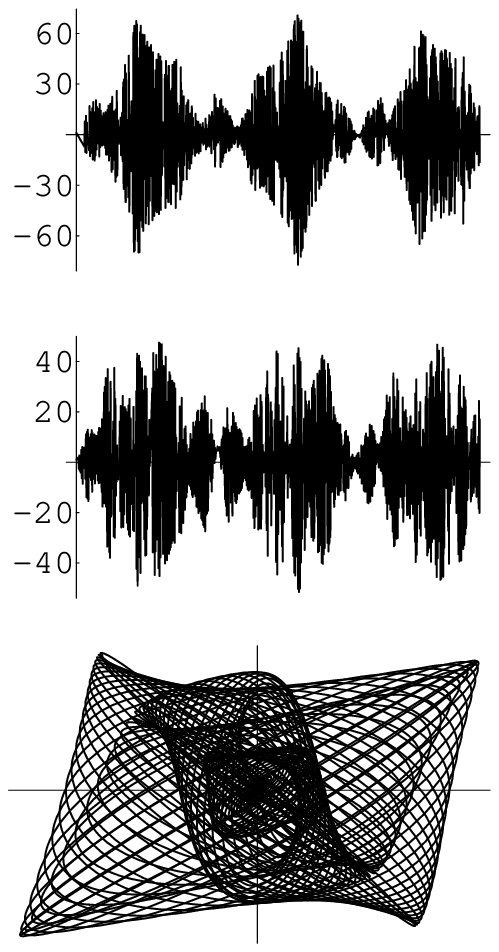}\hspace{0.5in}\includegraphics{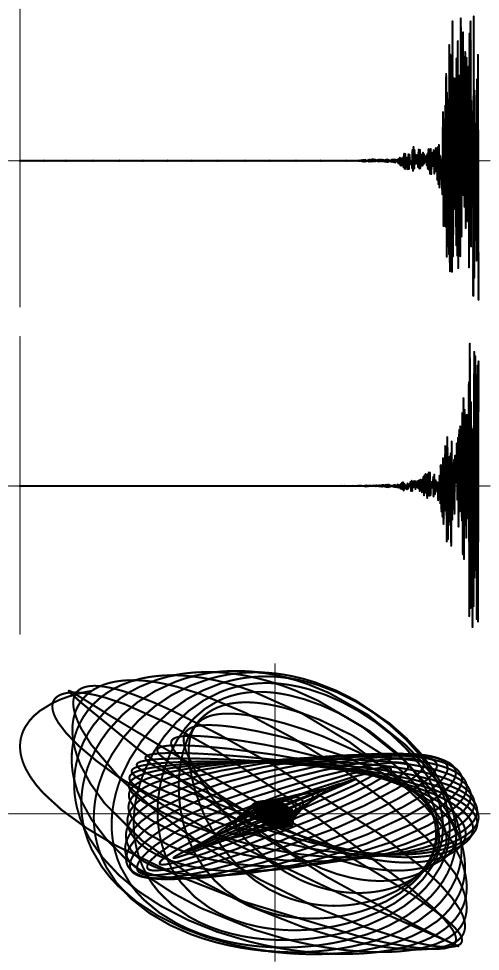}}
\end{center}
\caption{Long-time solutions for very close together points in
$\alpha,\beta$ parameter space. On the left is the solution for
$\alpha=435/800,\ \beta=425/800$; this point is in the stability
island indicated in Figure \ref{stabchart}. On the right is the
solution for $\alpha=436/800,\ \beta=425/800$; this solution is
outside the stability island. The integration times here are
$1600\pi$, and in both columns we present the $x$- and
$y$-solution versus time and the $x$-$y$ solution. The initial
conditions for both are
$x(0)=1,y(0)=1,\dot{x}(0)=0,\dot{y}(0)=0$.} \label{dpoint}
\end{figure*}

Based on these considerations, we will calculate the
multipliers of the monodromy matrix and find their norms; if the
norm is less than some cut-off value we say the system is stable,
otherwise unstable.

We focus on $(\alpha,\beta)=(0,1)\times(0,1)$ (as we may do due to
the scale invariance of the system, mentioned previously).
Splitting this interval into $n$ equal segments, we let
$\alpha=i/n$ and $\beta=j/n$ for $i,j=1,\ldots,n$. The fundamental
period of system \eqref{linsys} is given by \[ T=\frac{2\pi
n}{GCD[i,j]}
\] where $GCD$ is the greatest common divisor. We then integrate
the system \eqref{lininit} for $T$ units of time, evaluate
$\Phi(T)$, and calculate its eigenvalues.

The most important issue with regards accuracy is the long
integration time. If $i$ and $j$ do  not have a common divisor,
than $T=2\pi n$, which is large for $n$ large. We may halve this
integration time by taking advantage of the $t$-invariant nature
of the system (see above). If $A(t)$ is such that \eqref{tinv}
holds, then we may show that \[
\Phi(T)=G.\Phi(T/2)^{-1}.G.\Phi(T/2) \] since $G=G^{-1}$. Hence we
need only integrate over one half period to calculate $\Phi(T)$.
This improves the accuracy greatly, with one caveat: we must
invert $\Phi(T/2)$. However, since Louiville's theorem holds for
all $t$, this means the determinant of $\Phi(T/2)$ should be one,
and thus inverting may be done accurately.

In fact, we may use the condition $det[\Phi(t)]=1$ to monitor the
accuracy of  the integration, however this proves very costly with
regards integration time. Indeed, for some values of $\alpha$,
$\beta$ and $\vep$ we would expect the solution to be grossly
unstable. Monitoring the accuracy of such solutions is inefficient
as our concern is the stability-instability divide, rather than
how large these very large multipliers may be. As such, we are
satisfied that if $det[\Phi(t)]$ grows very large this is due to
the huge growth in very unstable combinations of $\alpha$ and
$\beta$; we will not waste effort with the accuracy of these
solutions.

In Figure \ref{stabchart} we give the stability diagram for system
\eqref{syseqn} with $\vep=0.1$.  Here, we set $n=800$ and thus
divide the $\alpha,\beta$ unit square parameter space into $800^2$
data points. For each we calculate the maximum multiplier norm,
with white regions being where the max norm is larger than 1.025,
black regions where the max norm is less than 1.025. We may make
some observations:

\begin{itemize}
  \item The most striking features are the large, linear bands of
  instability. These are aligned along lines of constant slope in
  the parameter space, in other words where $\alpha/\beta=$constant.
  These are clearly the dominant, low order resonances between
  $\alpha$ and $\beta$, and the lower the order of resonance the
  wider the band of instability. There are also higher order resonances,
  and these are of narrower width. This behavior is very much what
  we would expect, based on the classical Mathieu equation and its
  variations.
  \item There are two interesting and unexpected features. The first
  is the large `arc' of instability seen in the corners. This curved
  band of instability is not predicted by simple resonance between
  $\alpha$ and $\beta$.
  \item Second and perhaps more interesting are the small islands of stability
  which appear in the large $\alpha\approx\beta$ resonance band. We would
  fully expect the bands of instability to grow with $\vep$, and this
  we do observe. However, for small pockets of stability to develop
  when surrounding them is very large regions of instability is most unexpected.
  These islands are precarious: we show in Figure \ref{dpoint} the solution for
  two very nearby points in parameter space, one inside a small island
  of stability and one outside.
  \item The corresponding plot for other values of $\vep$ can be inferred from Figure
  \ref{stabchart} in the following way. Larger values of $\vep$ generate stability
  charts which `zoom in' towards the origin of Figure \ref{stabchart}, smaller values
  `zoom out'; this is due to the scale invariance of the system as described at the
  beginning of this section. For example, the stability chart for $\alpha,\beta$ on
  the unit square and $\vep=0.4$ is the same as the lower left quadrant of the same
  chart for $\vep=0.1$, since the parameter values $(\alpha,\beta,\vep)=(1/2,1/2,0.1)$
  have the same stability properties as $(1,1,0.4)$. This is why there is a cluster of
  unstable points near the origin: moving towards the origin is equivalent to increasing
  $\vep$; the resonance bands widen and overlap. As we {\it decrease} $\vep$,
  regions beyond the unit square in Figure \ref{stabchart} are drawn in; since the
  resonance bands narrow as we move away from the origin this is equivalent to saying they
  narrow as $\vep$ is decreased, as expected, and the chart is increasingly made up primarily
  of stable points. The exception is the $\alpha\approx\beta$ resonance band; it does not
  appear to taper as we move along it away from the origin. Some insight into this is gained in Section 4.
\end{itemize}

While the plots generated using numerical integration are
informative, they give us no feel for {\it why} certain bands of
instability develop or which are dominant. Nor do they give us an
analytic expression for the values at which the system bifurcates
from stable to unstable. In the next two sections, we use
analytical techniques which do just this, and what's more
recapture the true quasiperiodic nature of the original system.

\section{Transition curves using the method of harmonic balance}

In Floquet theory, $t$-invariant systems must have multipliers on
the unit circle or real line (with the rare exception of Figure
\ref{multipliers}(d)). As the system passes from stable to
unstable, the multipliers must pass through the point $(1,0)$ or
$(-1,0)$ in the complex plane. These multipliers represent
solutions which have the same period, or twice the period, of the
original system. This means solutions transiting from stable to
unstable have a known frequency and thus can be written as a
Fourier series. For this series to be a solution to the system
puts a constraint on the parameters involved, and this in turn
tells us the transition curves in parameter space.

In quasiperiodic systems, such as the one under consideration in
this paper, which we give again as \be
\ddot{x}+[\alpha^2+\vep\cos(\beta t)]x+\vep\cos(\alpha t)y =0
\nonumber \\ \ddot{y}+\vep \cos(\beta t)x+[\beta^2+\vep
\cos(\alpha t)]y=0, \label{seq} \ee Floquet theory does not hold.
However, based on the discussion above we can use the following ansatz:
solutions on transition curves from regions of stability to
instability in parameter space have the same quasiperiodic nature
as the system itself. This was justified and used by Rand et al in
\cite{rand1,rand4}, and while lacking in rigorous theoretical proof we find
this assertion works well in predicting the transition curves in
parameter space.

Since any quasiperiodic function can be written as an infinite
Fourier series in its frequencies (see, for example, Goldstein et
al \cite{goldstein}), we write the solutions to system \eqref{seq}
in the following form, \begin{align}
  \left(\begin{array}{c}
    x \\ y
  \end{array}\right)=\sum_{n=-\infty}^{\infty}\sum_{m=-\infty}^{\infty}
  \left(\begin{array}{c} A_{nm} \\ B_{nm} \end{array}\right)\exp\left[\frac{it}{2}(n\alpha+m\beta)\right], \label{fourser}
\end{align} where the 2 appears in the exponent to account for
solutions with twice the period of the system (corresponding to
the $-1$ multiplier). We substitute \eqref{fourser} into
\eqref{seq}, use the identity $\cos(\theta
t)=\tfrac{1}{2}(e^{i\theta t}+e^{-i\theta t})$, and shift indices
to cancel out the exponential terms. What remains are the
following two recurrence relations:
\begin{align} \left(\alpha^2-\frac{(n\alpha+m\beta)^2}{4}\right)A_{nm}\hspace{1.5in}\nonumber \\
+\frac{\vep}{2}\left(A_{n,m-2}+A_{n,m+2}+B_{n-2,m}+B_{n+2,m}\right)=0, \label{recurr1} \\
\left(\beta^2-\frac{(n\alpha+m\beta)^2}{4}\right)B_{nm}\hspace{1.5in}\nonumber \\
+\frac{\vep}{2}\left(A_{n,m-2}+A_{n,m+2}+B_{n-2,m}+B_{n+2,m}\right)=0,
\label{recurr2} \end{align} with $n,m=-\infty,\ldots,\infty$. This
is a set of linear, homogeneous equations in ${A_{nm},B_{nm}}$ of
the form \[ C \bmt{x}=\bmt{0}, \] and for non-trivial solutions to
exist we require $det(C)=0$; these are the well know infinite
Hill's determinants. In practice, we truncate at some order $N$ so
as $n,m=-N,\ldots,N$. Truncating at larger values of $N$ will
produce more accurate determinants, provided the determinant
converges. However, close inspection of the recurrence relations
(\ref{recurr1},\ref{recurr2}) reveals that this is an issue for
the system in question.

Writing the coefficient of $A_{nm}$ in \eqref{recurr1} as $\gamma_{nm}$ and dividing across gives the following: \begin{align}
  A_{nm}+\frac{\vep}{2\gamma_{nm}}\left(A_{n,m-2}+\ldots\right)=0
\end{align} with a similar expression for \eqref{recurr2}. In most systems, as $N\to\infty$ the `off-diagonal'
terms here grow progressively smaller, ensuring convergence. However, we see that if we let, for example,
$n=2,m=0$ then \[ \gamma_{20}=0. \] Thus the infinite determinant does {\it not} converge. This is a
problem particular to systems with a driving frequency that is equal to a natural frequency.

As a way around this problem we notice that a solution with period $2T$ also has period $4T$, and
indeed $6T,8T,\ldots$ and so on. Thus we may write \begin{align}
  \left(\begin{array}{c}
    x \\ y
  \end{array}\right)=\sum_{n=-\infty}^{\infty}\sum_{m=-\infty}^{\infty} \left(\begin{array}{c} A_{nm} \\
  B_{nm} \end{array}\right)\exp\left[\frac{it}{2M}(n\alpha+m\beta)\right],
\end{align} where $M=1,2,3,\ldots$. In this case we see \[ \gamma_{nm}=\left(\alpha^2-\frac{(n\alpha+m\beta)^2}{(2M)^2}\right) \]
and the singular term does not appear in the Hill determinant until $N=2M$. Thus we may work up to order
$N=2M-1$ and calculate the determinant. It will be a function of the parameters $(\alpha,\beta,\vep)$
and by setting it equal to zero we may therefore find analytic approximations to the transition curves in parameter space.

\begin{figure}
\begin{center}
\resizebox*{3in}{!}{\includegraphics{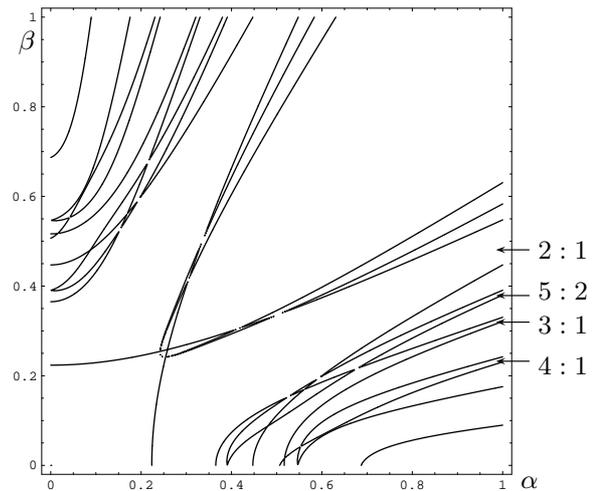}}
\end{center}
\caption{Some transition curves on the unit square for $\vep=0.1$ using the method of harmonic balance. The curves are where the determinant of the coefficient matrix $C$ vanishes.} \label{transcurves}
\end{figure}

For example, with $M=3$ we generate the matrix appropriate to 5th
order, find its determinant, and set it equal to zero. Factorizing
the determinant we see it contains terms like the following: \begin{eqnarray*}
&& (4\alpha^2-\beta^2+\vep)(4\alpha^2-\beta^2-\vep)=0,\nonumber \\
&& (9\alpha^4 - 82\alpha^2\beta^2 + 9\beta^4 - 64e^2)\times \\ && \hspace{0.5in}(9\alpha^4 - 82\alpha^2\beta^2 + 9\beta^4 -
16e^2)=0.\nonumber \end{eqnarray*} The first two clearly indicate a band of
instability with width $O(\vep)$ around the 2:1 resonance curve,
which grows to either side of the resonance as $\vep$ grows. The
second two terms indicate a band of width $O(\vep^2)$ around the
3:1 resonance curve, which grows {\it away} from the resonance as
$\vep$ grows. This implies that stability for resonant values of
$\alpha$ and $\beta$ is not as straightforward as we would
initially expect: for $\alpha/\beta=a/b$ with $a,b$ integers the
band of instability may move away from the resonance curve as
$\vep$ grows and thus commensurable values of $\alpha,\beta$ will
in fact be stable.

We give the zeroes of some of the main factors of the determinant
in Figure \ref{transcurves}. We can see in the figure some of the
main features we expect from our numerical simulations, however
some features are notable absent; foremost among these is the 1:1
resonance. Going to higher orders will capture more of the
transition curves, however the dimension of the problem quickly
becomes very large. If we truncate the series at $N=2M-1$ we must
find the determinant of a $2(2N+1)^2$ square matrix, which for the
example described above ($M=3$) is $242 \times 242$. Calculating
the determinant can be facilitated in the following way: these
large matrices are very sparse as is usual for recurrence
relations. We may separate out $2M$ submatrices, each describing
the $T,2T,\ldots,2MT$ solutions, as in \[ C=C_1 C_2 \ldots C_{2M},
\] and set the determinants of each to zero. However, the improvements in
accuracy from increasing $M$ do not justify the increasingly
cumbersome and lengthy matrices. Thus for this system it is
worthwhile to study the resonances of greatest import
individually; this we do in the next section.

\begin{center}
*******
\end{center}

Aside from these complications, we note that the method of harmonic balance can be used very quickly and directly to predict, in broad terms, the locations of the bands of instability. As mentioned previously, we may understand the mechanism of instability so: resonance curves in the $(\alpha,\beta)$ parameter space grow with $\vep$ into bands of instability due to resonance between the natural and driving frequencies. Predicting the location of said resonance curves will give approximately the regions of parameter space which will become unstable. As these resonance curves are the bands of instability in the limit of $\vep\to 0$, we may simply let $\vep=0$ in the recurrence relations given in (\ref{recurr1},\ref{recurr2}). This leads to the following definition of the resonance curves: \begin{align}
  \left(\alpha^2-\frac{(n\alpha+m\beta)^2}{4}\right)=\left(\beta^2-\frac{(n\alpha+m\beta)^2}{4}\right)=0,
\end{align} or \begin{align}
  \alpha=\frac{\pm m\beta}{2\pm n},\quad \beta=\frac{\pm n\alpha}{2\pm m}
\end{align} for $n,m=-\infty,\ldots,\infty$, unless $n(m)=\pm 2$ in which case $\beta(\alpha)=0$. These are precisely the straight lines in the $\alpha,\beta$ plane where the frequencies are in ratio.

\begin{figure}
\begin{center}
\resizebox*{3in}{!}{\includegraphics{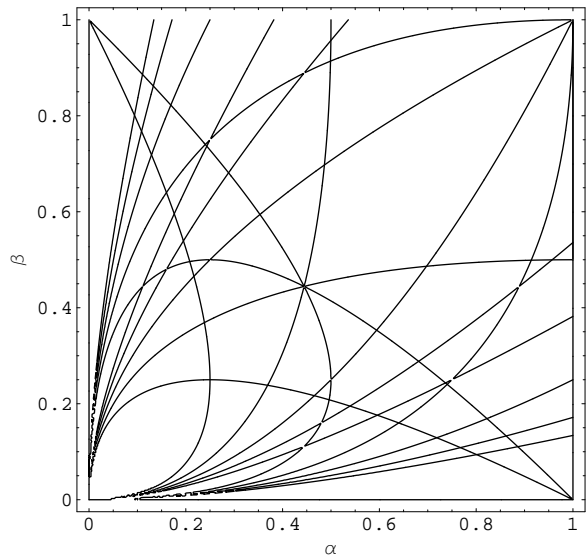}}
\end{center}
\caption{Resonance curves in the $\alpha,\beta$ plane for the
system \eqref{vareq}, with $|n|,|m|\leq 2$.} \label{rescurvevar}
\end{figure}

\begin{figure*}
\begin{center}
\resizebox*{4.7in}{!}{\includegraphics{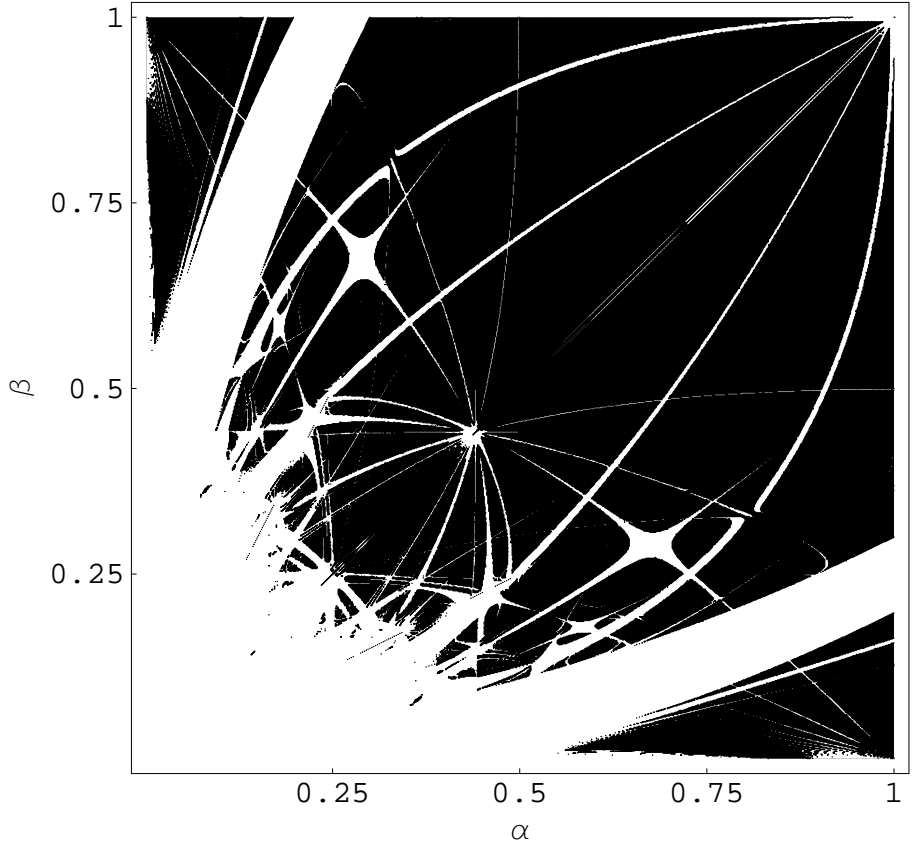}}
\end{center}
\caption{Stability diagram for system \eqref{vareq} with
$\vep=0.1$ and the $\alpha$, $\beta$ unit square divided into an
$800 \times 800$ grid. White represents points where the norm of the
maximum multiplier is greater than 1.025, and black where the norm
is less than 1.025. We see the strong correspondence between bands of instability in this plot and the resonance curves in Figure \ref{rescurvevar}.} \label{varstabchart}
\end{figure*}

As an illustration of how useful this technique can be we can examine the following problem, closely related to \eqref{seq}, \be \ddot{x}+[\alpha+\vep\cos(\beta
t)]x+\vep\cos(\alpha t)y =0 \nonumber \\ \ddot{y}+\vep \cos(\beta
t)x+[\beta+\vep \cos(\alpha t)]y=0. \label{vareq} \ee The time dependent terms here are of the same form as in \eqref{seq}, and thus we write solutions on the transition curves as in \eqref{fourser}. Now the recurrence relations, in the limit of $\vep\to 0$, yield the following resonance curves: \begin{align}
  \left(\alpha-\frac{(n\alpha+m\beta)^2}{4}\right)=\left(\beta-\frac{(n\alpha+m\beta)^2}{4}\right)=0,
  \label{reseqnvar}
\end{align} for $n,m=-\infty,\ldots,\infty$ (we see that in this case there is no issue with convergence of the infinite Hill's determinants, as the driving frequencies do not equal the natural frequencies). The resonance curves for the first few values of $n,m$ are given in Figure \ref{rescurvevar},
which can be compared with numerically generated stability charts
(using methods already discussed in Section 2) for non-zero $\vep$
as shown in Figure \ref{varstabchart}. We see that the resonance
curves as predicted by \eqref{reseqnvar} give a good approximation
to the location of the bands of instability in parameter space,
with very little effort.

In fact, we may generalize this technique to arrive at the following proposition:\\

\noindent{\bf Proposition:} For the system \be \frac{d^2X}{dt^2}+\left(P_0+\vep P_1(t)\right)X=0 \label{prop1} \ee
with $X\in\mathbb{R}^{n}$, $P_0=diag(\lambda_1^2,\lambda_2^2,\ldots,\lambda_n^2)$ and $P_1(t)$
of the form \bes P_1(t)= \left(
                   \begin{array}{ccc}
                     \cos(\omega_1 t) & \ldots & \cos(\omega_n t) \\
                     \vdots & \vdots & \vdots \\
                     \cos(\omega_1 t) & \ldots & \cos(\omega_n t) \\
                   \end{array}
                 \right),
 \ees  bands of instability will grow with $\vep$ from resonance curves given by \[ \lambda_j^2-\frac{(\bmt{n}.\bmt{\omega})^2}{4}=0,\quad j=1,\ldots,n, \]
where $\bmt{n}$ is a vector of length $n$ whose elements are all integers in the range $-\infty,\ldots,\infty$, that is $\bmt{n}\in\mathbb{Z}^n$.\\

\noindent{\bf Proof:} On transition curves, the solution has the same quasiperiodic structure as
$P_1(t)$ and so we can write it as a Fourier series of the form
\[ X(t)=\sum_{\bmt{n}} \bmt{A}_{\bmt{n}} \exp\left[\frac{it}{2} (\bmt{n}.\bmt{\omega}) \right] \]
where $\bmt{A}_{\bmt{n}}$ is an $n$-dimensional vector of Fourier coefficients. Subbing this into \eqref{prop1} we find

\begin{align*}
  &-\frac{(\bmt{n}.\bmt{\omega})^2}{4}\sum_{\bmt{n}} \bmt{A}_{\bmt{n}} \exp\left[\frac{it}{2} (\bmt{n}.\bmt{\omega}) \right] \\ &+(P_0+\vep P_1(t)) \sum_{\bmt{n}} \bmt{A}_{\bmt{n}} \exp\left[\frac{it}{2} (\bmt{n}.\bmt{\omega}) \right] =0.
\end{align*}

The term in $\vep$ in this equation is complicated and would
involve shifting indices were we looking for the infinite Hill's
determinant to determine the transition curves; however we seek
only the {\it resonance curves} which are found by letting $\vep=0$.
Canceling the exponentials leaves \[
\sum_{\bmt{n}}\left(-\frac{(\bmt{n}.\bmt{\omega})}{4}\bmt{A}_{\bmt{n}}+P_0\bmt{A}_{\bmt{n}}
\right)=0 \] and since $P_0=diag(\lambda_i^2)$ for $i=1,\ldots,n$
this leaves us with the $n$ equations \[
\lambda_j^2-\frac{(\bmt{n}.\bmt{\omega})^2}{4}=0, \quad \bmt{n}\in\mathbb{Z}^n. \] Now we may simply plot these curves for some elements of $\bmt{n}$ as in Figure \ref{rescurvevar}.

We note this method predicts the resonance points for the Mathieu equation \eqref{math1d} to be $\delta-n^2/4=0$, and the resonance curves for the quasiperiodic Mathieu equation \eqref{quasi1d} to be $\delta-(n+m\omega)^2/4=0$ (as in \cite{rand1}), with $n,m\in\mathbb{Z}$.

\section{The multiple scales method near the 1:1 resonance}

Based on the previous section we see that the method of harmonic
balance is unable to describe behavior near the $1:1$ resonance,
that is close to $\alpha=\beta$. We use instead the method of
multiple scales, which is a perturbative technique that allows us
to derive analytic representations near resonant values of
$\alpha$ and $\beta$.

The typical way to proceed is to expand $\alpha,\beta, x$ and $y$
in powers of $\vep$, then define different time scales and
separate out the behaviour on slower time scales. We could let \[
\alpha=\alpha_0+\vep\alpha_1+\vep^2\alpha_2,\quad
\beta=\alpha_0+\vep\beta_1+\vep^2\beta_2 \] so $\alpha$ and
$\beta$ only differ at first order. Then, the differential
equations would contain terms involving $\alpha_0 t+\alpha_1 \vep
t+\alpha_2 \vep^2 t$ which suggests we define three time scales,
\[ t, \quad \tau=\vep t,\quad T=\vep^2 t \] ($T$ here is not to be confused with the fundamental period $T$ from Section 2). Note the $\vep^2$
time scale is necessary as there are no secular terms at first
order. Now we may expand $x$ and $y$ in powers of $\vep$ where
each term is a function of these three times, i.e.
\begin{align} x=x_0(t,\tau,T)+\vep x_1(t,\tau,T)+\vep^2 x_2(t,\tau,T)+\ldots \label{xexpans} \end{align} and
similarly for $y$. Subbing into the differential equations and
collecting powers of  $\vep$ we find a set of equations with
inhomogeneous terms involving partial derivatives of lower order
solutions. Removing secular terms gives differential equations in
the slow-time dependence of low order terms. However, we find that
the crucial system of equations only describes flow on the `very'
slow time $T$, and is independent of $\alpha_1$ and $\beta_1$ (see the
Appendix for details).

Therefore, it is sufficient to define \begin{align} \alpha=\alpha,\quad
\beta=\alpha+\vep^2\beta_2. \label{albems} \end{align} Now there is only one
slow time, $T=\vep^2 t$, however we are still expanding $x$ and
$y$ in powers of $\vep$ (rather than $\vep^2$). To second order we
find the following system of equations:

\begin{align*}
  \vep^0: & \frac{\partial^2x_0}{\partial t^2}+\alpha^2 x_0=0, \\
          & \frac{\partial^2y_0}{\partial t^2}+\alpha^2 y_0=0, \\
  \vep^1: & \frac{\partial^2x_1}{\partial t^2}+\alpha^2 x_1=-x_0
  \cos(\alpha t+\beta_2T)-y_0 \cos(\alpha t), \\
          & \frac{\partial^2y_1}{\partial t^2}+\alpha^2 y_1=-y_0
  \cos(\alpha t+\beta_2T)-y_0 \cos(\alpha t), \\
  \vep^2: & \frac{\partial^2x_2}{\partial t^2}+\alpha^2 x_2=-2 \frac{\partial^2x_0}{\partial t\partial T}-x_1
  \cos(\alpha t+\beta_2T)\\ &\hspace{1in}-y_1 \cos(\alpha t), \\
          & \frac{\partial^2y_2}{\partial t^2}+\alpha^2 y_2=-2 \frac{\partial^2y_0}{\partial t\partial T}-y_1
  \cos(\alpha t+\beta_2T)\\ &\hspace{1in}-y_1 \cos(\alpha t)-2\alpha\beta_2 y_0, \\
\end{align*} where $x_0=x_0(t,T)$ and so on. Solving one at a time,
 we find at zero order: \begin{subequations} \label{zer12} \be x_0=A(T)\cos(\alpha t)+B(T)\sin(\alpha
t),\\ y_0=C(T)\cos(\alpha t)+D(T)\sin(\alpha t). \ee \end{subequations} There are no
secular terms on the right hand side at first order so we find \[
x_1=a(T)\cos(\alpha t)+b(T)\sin(\alpha
t)+\frac{A}{6\alpha^2}\cos(2\alpha t+\beta_2T)+\ldots \] and a
similar expression for $y_1$. Now when we look at second order, we
see there are secular terms in the right hand side (that is, terms
involving $\cos(\alpha t)$ and $\sin(\alpha t)$). Setting these
terms to zero results in four linear differential equations in the slow
time dependence of $x_0$ and $y_0$, that is the four functions
$A(T),B(T),C(T)$ and $D(T)$ defined in \eqref{zer12}. Letting $Y=(A,B,C,D)$ we may write
these four equations as $dY/dT=Q(T)Y$ where

\begin{strip}
\[
Q=\left(%
\begin{array}{cccc}
  -\sin(\beta_2T)+3\sin(2\beta_2T) & -2+\cos(\beta_2T)+3\cos(2\beta_2T) & 7\sin(\beta_2T) & 1+\cos(\beta_2T) \\
  2+5\cos(\beta_2T)+3\cos(2\beta_2T) & -7\sin(\beta_2T)-3\sin(2\beta_2T) & 5+5\cos(\beta_2T) & \sin(\beta_2T) \\
 -\sin(\beta_2T)+3\sin(2\beta_2T) & -2+\cos(\beta_2T)+3\cos(2\beta_2T) & 7\sin(\beta_2T) & 1+24\alpha^3\beta_2+\cos(\beta_2T) \\
  2+5\cos(\beta_2T)+3\cos(2\beta_2T) & -7\sin(\beta_2T)-3\sin(2\beta_2T) & 5-24\alpha^3\beta_2+5\cos(\beta_2T) & \sin(\beta_2T) \\
\end{array}%
\right) \]
\end{strip}

The stability of this system will determine the
stability of the original system.

\begin{figure}
\begin{center}
\resizebox*{3in}{!}{\includegraphics{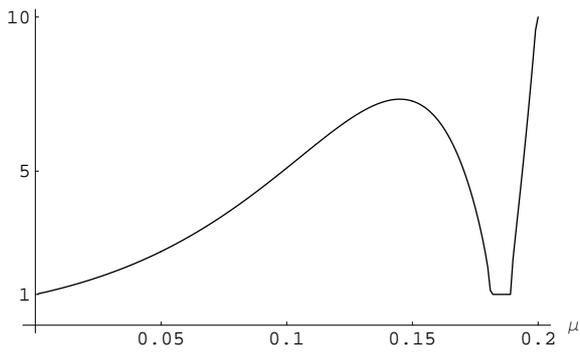}}
\end{center}
\caption{The largest multiplier of system \eqref{musys} for $\mu$ small.} \label{scalstab}
\end{figure}

This system is still Louiville ($Tr(Q)=0$) and $t$-invariant (as defined in \eqref{tinv} with $G=diag(1,-1,1,-1)$). A nice simplification happens if we define a new time coordinate $t^*=\beta_2 T$, and we may write the system in the form \begin{align}
  \frac{dX}{dt^*}=&\left[ A+\mu (A_0+ A_1 \cos(t^*)+ A_2 \cos(2t^*)\right. \nonumber \\ &\left. \hspace{0.5in} +B_1 \sin(t^*)+B_2 \sin(2t^*))\right]X \label{musys}
\end{align} where $\mu=1/(24\alpha^3\beta_2)$ and the constant coefficient matrices are given by \begin{align*} \setlength\arraycolsep{0.05in}
  A=\left(
      \begin{array}{cccc}
        0 & 0 & 0 & 0 \\
        0 & 0 & 0 & 0 \\
        0 & 0 & 0 & 1 \\
        0 & 0 & -1 & 0 \\
      \end{array}
    \right), \quad A_0=\left(
      \begin{array}{cccc}
        0 & -2 & 0 & 1 \\
        2 & 0 & 5 & 0 \\
        0 & -2 & 0 & 1 \\
        2 & 0 & 5 & 0 \\
      \end{array}
    \right),
\end{align*} with $A_1,A_2,B_1,B_2$ in the Appendix.
This is a one parameter $2\pi$-periodic system, with the $\mu$ parameter multiplying terms very much like a Fourier series. It is straightforward to analyze the stability of this system using the methods already outlined, and we show the maximum multiplier in Figure \ref{scalstab}.

\begin{figure}
\begin{center}
\resizebox*{3in}{!}{\includegraphics{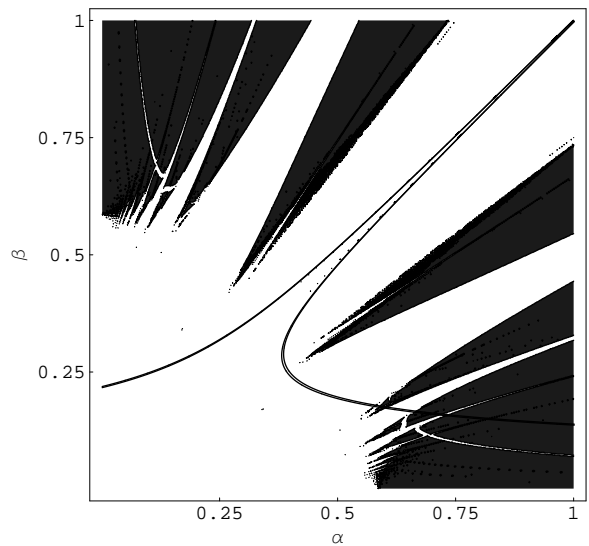}}
\end{center}
\caption{The multiple scales prediction for a band of stability overlaid onto the stability chart generated numerically in Section 2.} \label{scalover}
\end{figure}

The appearance of a small region of stability can be understood in the following way. In system \eqref{musys}, when $\mu=0$ the system has natural frequency 1. As we switch on $\mu$, the driving frequency is also 1 and resonance causes the system to be unstable. However, as $\mu$ grows, the natural frequencies are now the eigenvalues of $A+\mu A_0$, which become sufficiently distinct from the driving frequency for stability to resume. This is transient however and as the `perturbation' $\mu$ grows the system becomes unstable again. This window of stability occurs for \[ (\mu_-=)\ 0.18106\lesssim\mu\lesssim 0.189\ (=\mu_+), \] which tells us the values of $\beta_2$ where stability occurs and this in turn gives us the following band of stability in $(\alpha,\beta)$ parameter space \[ \alpha\pm\vep^2\left(\frac{1}{24\alpha^3\mu_+}\right)<\beta<\alpha\pm\vep^2\left(\frac{1}{24\alpha^3\mu_-}\right) \] This result is in excellent agreement with the numerical results generated earlier, as shown in Figure \ref{scalover}, at least for values of $\alpha,\beta$ over 0.5. As discussed previously, this is because small values of $\alpha,\beta$ are equivalent to large values of $\vep$ and the multiple scales method would not be expected to work well in this regime. Note that the symmetric nature of the system is recovered by the expansion \begin{align} \beta=\beta,\quad
\alpha=\beta+\vep^2\alpha_2, \end{align} analogous to \eqref{albems}, which results in the same stability region as in Figure \ref{scalover} reflected about the line $\alpha=\beta$.

\section{Conclusions}

In this paper we have begun to analyze the rich dynamical structure of two dimensional systems with quasiperiodic coefficients. Through a combination of numerical and perturbative techniques we are able to predict and determine the transition curves in parameter space which demarcate the regions leading to stable and unstable solutions. In particular, we find that using the method of harmonic balance we can quickly approximate the resonance curves in parameter space from which bands of instability emanate, with the caveat that the band of instability may spread out on either side of the resonance curve or move to its side. The multiple scales method has proven to be particularly informative, if a little cumbersome, for examining the vicinity of specific resonances.

There are three obvious issues which need to be understood further with regards the particular system under consideration in this work. The first is the `arcs' of instability which are commented on in \S 2 and can be seen in Figure \ref{stabchart}. These arcs do not arise in the harmonic balance analysis, nor do they correspond to a linear relationship between the frequencies and thus the method of multiple scales would be inappropriate. The exact mechanism which leads to the development of these arcs is not understood, and is worthy of further investigation.

The second issue is the possible resonances between first order terms in the expansion of $\alpha,\beta$ as described in the Appendix. Perhaps these could explain why the multiple scales method predicts a continuous {\it band} of stability near the 1:1 resonance, whereas numerically we observe a broken series of islands of stability (see Figure \ref{stabchart} and \ref{scalover}). Finally, all bands of instability visible in Figure \ref{stabchart} taper as we increase $\alpha$ and $\beta$, as predicted by the harmonic balance method, with the exception of the large region of instability surrounding the 1:1 resonance. While the multiple scales analysis could be interpreted as telling us that the true band of instability around $\alpha=\beta$ is in fact very narrow and exists between the predicted bands of {\it stability}, as shown in Figure \ref{scalover}, it still leaves an open question as to the origin of this large region of instability.

There are numerous extensions to the system in question, however it would be wise to keep the number of free parameters small to enable an informative analysis. For example, a general form could consist of two natural frequencies and four driving frequencies. While this would undoubtedly lead to exceedingly rich stability properties, six free parameters would be difficult to present in a clear fashion. A `detuning' parameter in the driving frequencies would instead lead to three parameters and would provide an interesting extension. There are of course other extensions involving nonlinear terms, more dimensions and so on. Finally, there are numerous applications involving quasiperiodic oscillations about fixed points whose analysis would benefit from the present work.

\section{Appendix}

As is typical in multiple scales analysis, the equations at high
orders are quite lengthy and involve very many terms; as such we
will only give an overview to justify ignoring first order terms
in $\alpha,\beta$ to arrive at \eqref{albems}.

Let us write  \[ \alpha=\alpha_0+\vep\alpha_1+\vep^2\alpha_2,\quad
\beta=\alpha_0+\vep\beta_1+\vep^2\beta_2 \] and expand $x,y$ in
the three time scales $t,\tau$ and $T$, as in \eqref{xexpans}. The
second time derivative becomes (up to second order in $\vep$) \[
\frac{d^2}{dt^2}=\frac{\partial^2 }{\partial
t^2}+\vep\left(2\frac{\partial^2 }{\partial t \partial \tau }
\right)+\vep^2\left(2\frac{\partial^2 }{\partial t\partial T
}+\frac{\partial^2 }{\partial \tau^2} \right), \] leading to the
following system of equations in $x$:

\begin{align*}
  \vep^0:& \frac{\partial^2 x_0}{\partial t^2}+\alpha_0^2 x_0=0 \\
  \vep^1:& \frac{\partial^2 x_1}{\partial t^2}+2\frac{\partial^2 x_0}{\partial t\partial
  \tau}+\alpha_0^2 x_1+2\alpha_0\alpha_1 x_0 \\ & +x_0 \cos(\alpha_0
  t+\beta_1\tau+\beta_2T) \\ & +y_0\cos(\alpha_0t+\alpha_1\tau+\alpha_2T)=0
  \\
  \vep^2:& \frac{\partial^2 x_2}{\partial t^2}+2\frac{\partial^2 x_1}{\partial
  t\partial\tau}+2\frac{\partial^2 x_0}{\partial t\partial T}+\frac{\partial^2 x_0}{\partial
  t^2} \\ & +\alpha_0 x_2 +2\alpha_0\alpha_1
  x_1+(\alpha_1^2+2\alpha_0\alpha_2)x_0 \\ & +x_1\cos(\alpha_0
  t+\beta_1\tau+\beta_2T) \\ & +y_1\cos(\alpha_0t+\alpha_1\tau+\alpha_2T)=0
\end{align*}

For the corresponding system in $y$, simply swap $x$ with $y$, and
$\alpha_1,\alpha_2$ with $\beta_1,\beta_2$. We examine the
equations one order at a time. At zero order, the solutions are
\begin{align*}
  x_0=a(\tau,T)\cos(\alpha_0 t)+b(\tau,T)\sin(\alpha_0 t), \\
y_0=c(\tau,T)\cos(\alpha_0 t)+d(\tau,T)\sin(\alpha_0 t).
\end{align*} At first order, which we write as \[ \frac{\partial^2x_1}{\partial t^2}+\alpha_0 x_1=\textrm{inhomogeneous terms}, \] the right hand side
includes
\begin{align*}
  2\alpha_0\sin(\alpha_0 t)\left(\frac{\partial a}{\partial \tau}-\alpha_1 b \right)
  +  2\alpha_0\cos(\alpha_0 t)\left(-\frac{\partial b}{\partial \tau}-\alpha_1 a \right)
\end{align*} in the $x$-equation, and similarly for $c,d$ in the
$y$-equation. Setting these secular terms equal to zero gives us
the $\tau$-dependence at zero order, i.e. \begin{align*}
a(\tau,T)=A(T) \cos(\alpha_1 \tau)+A'(T)\sin(\alpha_1\tau), \\
b(\tau,T)=B(T) \cos(\alpha_1 \tau)+B'(T)\sin(\alpha_1\tau),
\end{align*} and similarly for $c,d$. The solution at first order will then involve the
homogeneous solution \begin{subequations} \label{firstordsoln} \begin{align} x_1=P(\tau,T)\cos(\alpha_0
t)+Q(\tau,T)\sin(\alpha_0 t), \\ y_1=R(\tau,T)\cos(\alpha_0
t)+S(\tau,T)\sin(\alpha_0 t), \end{align} \end{subequations} and the inhomogeneous
solution, made up of $\sin$ and $\cos$ terms with arguments \[
(\alpha_0\pm\alpha_0)t+(\alpha_1/\beta_1)\tau+(\alpha_2/\beta_2)T
\] which we will not list in full (by $\alpha_i/\beta_i$ we mean linear combinations of the form $k_1\alpha_i+k_2\beta_i$ with $k\in\mathbb{Z}$ and $|k_1|+|k_2|\leq 3$).

At second order, which again we write as \[
\frac{\partial^2x_2}{\partial t^2}+\alpha_0
x_2=\textrm{inhomogeneous terms}, \] we use the zero and first
order solutions to calculate the right hand side and look for
secular terms, that is terms multiplying $\cos(\alpha_0
t),\sin(\alpha_0 t)$. The full expression for the right hand side
is very lengthy and complicated, and involves multiplying many
$\cos$ and $\sin$ terms of different arguments together, however
we are only interested in those with $\alpha_0$ in the
$t$-dependence and so can ignore those containing $2\alpha_0$ or
$3\alpha_0$. This gives us four equations: the coefficients of
$\cos(\alpha_0 t),\sin(\alpha_0 t)$ in the $x$-equation, and the
same for the $y$-equation, which must be set equal to zero to
prevent resonance in $x_2,y_2$. These four equations are of the
form \begin{subequations} \label{secord} \begin{align} -2\alpha_0 \frac{\partial Q}{\partial
\tau}-2\alpha_0\alpha_1P+f_1(t,\tau,T)=0, \\ 2\alpha_0
\frac{\partial P}{\partial
\tau}-2\alpha_0\alpha_1Q+f_2(t,\tau,T)=0,
\end{align} \end{subequations} and again in $R,S$ and $f_3,f_4$ with $\alpha_1$ replaced by $\beta_1$. These define the
$\tau$-dependence of $P,Q,R,S$ (from the first order solutions
$x_1,y_1$ in \eqref{firstordsoln}), and we see these will be oscillatory if we remove the
$\cos(\alpha_1\tau)$,$\sin(\alpha_1\tau)$ terms in
$f_1,f_2$, and the $\cos(\beta_1\tau)$,$\sin(\beta_1\tau)$ terms in
$f_3,f_4$. This leads to eight equations involving
the eight functions describing the $T$-dependence of the zero
order solutions, that is \[ A,A',B,B',C,C',D,D' \] and their
derivatives with respect to $T$. We may write this as a first
order system
$dZ/dT=\Lambda(T)Z$, where $\Lambda$ has block form \[ \Lambda(T)=\left(%
\begin{array}{cc}
  \Lambda_1 & \Lambda_2(T) \\
  \Lambda_3(T) & \Lambda_4 \\
\end{array}%
\right), \] with $\Lambda$ being Louiville (vanishing trace) and
$t$-invariant, as defined in \eqref{tinv} with
\[ G=diag(1,-1,-1,1,1,-1,-1,1).\] The stability of our original
system now hangs on the stability of this one, and crucially the
parameters $\alpha_1,\beta_1$ do not appear in the matrix
$\Lambda$. Thus the difference between $\alpha$ and $\beta$ at
first order is irrelevant and we may set $\alpha_1=\beta_1=0$.

We note that there are a number of possible resonances in the first order terms in $\alpha,\beta$ which we have ignored. The $f_i$ functions described above contain terms like, for example, \begin{equation}
  \cos(2\alpha_1 \tau+\beta_1\tau+\alpha_2T-\beta_2T).
\end{equation} In general, this will not be secular in equations \eqref{secord}. If however it was the case that $\beta_1=-\alpha_1$ then this {\it would} be secular. There are six combinations which lead to additional secular terms, they are \[ \alpha_1=\pm\beta_1,\quad \alpha_1=\pm 3\beta_1,\quad 3\alpha_1=\pm\beta_1. \] Each would need to be considered in turn, and would lead to a different coefficient matrix $\Lambda$ given above. We will not pursue this issue further in the present work.

Finally, we give the coefficient matrices from \eqref{musys},
\begin{align*} \setlength\arraycolsep{0.06in} A_1=\left(
      \begin{array}{cccc}
        0 & 1 & 0 & 1 \\
        5 & 0 & 5 & 0 \\
        0 & 1 & 0 & 1 \\
        5 & 0 & 5 & 0 \\
      \end{array}
    \right), \quad A_2=\left(
      \begin{array}{cccc}
        0 & 3 & 0 & 0 \\
        3 & 0 & 0 & 0 \\
        0 & 3 & 0 & 0 \\
        3 & 0 & 0 & 0 \\
      \end{array}
    \right), \\ \setlength\arraycolsep{0.04in} B_1=\left(
      \begin{array}{cccc}
        -1 & 0 & 7 & 0 \\
        0 & -7 & 0 & 1 \\
        -1 & 0 & 7 & 0 \\
        0 & -7 & 0 & 1 \\
      \end{array}
    \right),  B_2=\left(
      \begin{array}{cccc}
        3 & 0 & 0 & 0 \\
        0 & -3 & 0 & 0 \\
        3 & 0 & 0 & 0 \\
        0 & -3 & 0 & 0 \\
      \end{array}
    \right).
\end{align*}


\bibliographystyle{plain}
\bibliography{mathieubib}   

\end{document}